%% file: main.tex
\documentclass[12pt]{amsart}
\usepackage{enumerate,amssymb,amsfonts,latexsym,psfrag,psfig}
\usepackage{amscd,amsmath}
\usepackage[dvips]{graphicx}

\input{macros.tex}
\begin{document}

\bigskip\bigskip

\title[renormalizationof area-preserving maps]{
Rigidity for infinitely renormalizable area-preserving maps\\}

\author{D. Gaidashev, T. Johnson , M. Martens}

\address {Department of Mathematics, Uppsala University, Uppsala, Sweden, gaidashev@math.uu.se}

\address{Fraunhofer-Chalmers Research Centre for Industrial Mathematics
Chalmers University of Technology, 
SE-412 88 Gothenburg, Sweden,
tomas.johnson@fcc.chalmers.se}

\address{Department of Mathematics, Stonybrook University, Stony Brook, NY11794, USA, marco@math.sunysb.edu}

\date{\today}

\begin{abstract} 
The period doubling Cantor sets of strongly dissipative H\'enon-like maps with different average Jacobian are not smoothly conjugated. The Jacobian Rigidity Conjecture says that the period doubling Cantor sets of two-dimensional H\'enon-like  maps with the same average Jacobian are smoothly conjugated. This conjecture is true for average Jacobian zero, e.g. the one-dimensional case.  The other extreme case is when the maps preserve area, e.g. the average Jacobian is one. Indeed, the period doubling Cantor set of area-preserving maps in  the universality class of the Eckmann-Koch-Wittwer renormalization fixed point
are smoothly conjugated.
\end{abstract}

\maketitle

\setcounter{tocdepth}{1}
\tableofcontents

\input{intro}

\input{prelim}
\input{cantorset}
\input{rigidity}
\input{appendix}

\input{references}

\end{document}

%% file: macros.tex
\newtheorem{thm}{Theorem}[section]

\newtheorem{lem}[thm]{Lemma}
\newtheorem{prop}[thm]{Proposition}
\theoremstyle{remark}
\newtheorem{rem}{Remark}[section]

\theoremstyle{definition}
\newtheorem{defn}{Definition}[section]

\numberwithin{equation}{section}
\numberwithin{figure}{section}

 \newcommand{\pichere}[2]
{\begin{center}\includegraphics[width=#1\textwidth]{#2}\end{center}}

\font\nt=cmr7

\def\note#1
{\marginpar
{\nt $\leftarrow$
\par
\hfuzz=20pt \hbadness=9000 \hyphenpenalty=-100 \exhyphenpenalty=-100
\pretolerance=-1 \tolerance=9999 \doublehyphendemerits=-100000
\finalhyphendemerits=-100000 \baselineskip=6pt
#1}\hfuzz=1pt}

\newcommand{\ra}{\rightarrow}

\def\ssk{\smallskip}

\newcommand{\diam}{\operatorname{diam}}
\newcommand{\dist}{\operatorname{dist}}

\newcommand{\inter}{\operatorname{int}}
\renewcommand{\mod}{\operatorname{mod}}

\newcommand{\id}{\operatorname{id}}

\newcommand{\transverse}
 {\kern .7em\makebox[0pt][c]{\raisebox{.2ex}{$|$}}\kern -.6em\cap}

\newcommand{\tangent}
 {\kern .7em\makebox[0pt][c]{\raisebox{.77ex}{$--$}}\kern -.6em\cap}

\newcommand{\OO}{{\mathcal O}}

\newcommand{\R}{{\mathbb R}}

\newcommand{\Z}{{\mathbb Z}}

\def\BPhi{{\boldsymbol{\BPhi}}}
\def\B0{{\mathbf{0}}}

\newcommand{\Dom}{\operatorname{Dom}}

\catcode`\@=12

\def\Empty{}
\newcommand\oplabel[1]{
  \def\OpArg{#1} \ifx \OpArg\Empty {} \else
  	\label{#1}
  \fi}
		
\newcommand{\comm}[1]{}
\newcommand{\comment}[1]{}

%% file: intro.tex
\section{Introduction}\label{intro}

The simplest dynamics beyond periodic behavior can be found in systems at the accumulation of period doubling. These systems have a period doubling Cantor set whose dynamics can be very well understood in terms of the periodic orbits which accumulate on this set. The consecutive approximating periodic orbits are related by so-called period doubling bifurcations. This simplest non-period behavior is often observed in models and in nature. In particular, the transition to chaos in strongly dissipative systems often occurs with such a period doubling Cantor set.

The topological properties of these period doubling Cantor sets are simple and very well understood.   In the mid-1970's 
Feigenbaum~\cite{F1,F2} and, independently,  Coullet and
Tresser \cite{CT,TC}, discovered universal geometrical properties
of these Cantor sets at the transition to chaotic behavior in one-dimensional dynamics. Coullet and Tresser conjectured 
that the universal geometry at transition to chaos in one-dimensional dynamics will also be observed in higher dimensional systems. 
This conjecture has been confirmed by many numerical and physical 
experiments.\footnote{This conjecture should be taken with caution as not every transition to chaos is related to accumulation of period doubling.} 

 \bigskip
 
Feigenbaum-Coullet-Tresser introduced renormalization in dynamics to explain the observed geometrical universality. The renormalization operator acts as a microscope. By repeatedly applying, one can describe the small scale geometrical structure
by renormalizations, e.g. systems associated to the smaller scales.
The corresponding renormalization operator has a unique hyperbolic fixed point. The dynamics of the renormalization fixed point, which is itself a one-dimensional system, and the behavior of the renormalization operator around this fixed point determine the  asymptotic small scale properties. This explains the observed universality.

The renormalization technique has been generalized to many types of dynamics. However, a rigorous study of universality  has been surprisingly difficult and technically sophisticated and so far has only been thoroughly carried out in the case of one-dimensional maps, on the interval or the circle, see~\cite{AL, FMP,He, L,Ma,McM,MS, S,VSK,Y}. There is, however, at present  no deep understanding of universality in conservative systems, other than in the  case of the universality for systems ``near integrability'' \cite{AK,AKW,EKW2, K1,K2,K3,G1,Ko,Kh}. 

Convergence of renormalization implies universal geometry at asymptotic small scale around a certain point. A stronger property is rigidity: there is universal geometry at asymptotic small scale around all points. All period doubling Cantor sets are topologically equivalent. Between two of them there will be a homeomorphism conjugating
the dynamics.  A priori there is no reason to belief that these homeomorphisms have some smoothness. Rigidity means that there are conjugating homeomorphism which are smooth. The topology determines the small scale geometry.

\bigskip

Indeed, the period doubling Cantor sets in one-dimensional dynamics are rigid. 
This is surprising. Especially because the Cantor sets have very rich geometry: there are essentially no two points with the same asymptotic small scale geometry, see \cite{BMT}.
Many numerical and physical experiments show that exactly the same universal geometry from one-dimensional dynamics occurs also in some dissipative higher dimensional systems. Surprisingly, the rigidity phenomenon is more delicate in higher dimensions.

Strongly dissipative two-dimensional H\'enon-like maps can be thought of as two dimensional perturbations of one-dimensional systems. In \cite{CEK1} and \cite{CLM}  two renormalization schemes were developed for strongly dissipative H\'enon-like maps at the accumulation of period doubling. This explained the observed one-dimensional
universal geometry present in these H\'enon-like maps. Surprisingly, the period doubling H\'enon-like Cantor sets are not smoothly conjugated to their one-dimensional counter part, see \cite{CLM, LM1}.

The average Jacobian, e.g. the average rate of dissipation, plays a role. The Cantor set of systems with different average Jacobian 
are not smoothly conjugated, see \cite{CLM, LM1}. Nevertheless, the  
geometry of the one-dimensional period doubling Cantor set is still present. The conjugations between the Cantor sets of strongly dissipative H\'enon-like maps is almost everywhere, with respect to the natural measure on the Cantor set,  smooth. This phenomenon is called {\it Probabilistic Rigidity}, see \cite{LM1, LM3}. Small scale geometry has a probabilistic nature in higher dimensions.

As we discussed before, different average Jacobians are an obstruction to rigidity. However, this is not contradicting the rigidity paradigm {\it topology determines geometry}. Namely, in \cite{LM2} it is shown that 
maps with different average Jacobian are not topologically equivalent.
The average Jacobian is a topological invariant. It describes 
 topological aspects of  the heteroclinic web, e.g. the network of stable and unstable manifolds of the period orbits. This leads to the

\bigskip

\noindent
{\bf Jacobian Rigidity Conjecture.} {\it The period doubling Cantor sets of two-dimensional H\'enon-like maps with the same average Jacobian are smoothly conjugated.}

\bigskip

This conjecture is true for average Jacobian zero, e.g. the one-dimensional case.  The other extreme case when the maps preserve area, e.g. the average Jacobian is one,  is described in Theorem \ref{rigthm}. It says,

\bigskip

\noindent
{\bf Rigidity for Area-preserving Maps.} {\it The period doubling Cantor set of area-preserving maps in  the universality class of the Eckman-Koch-Wittwer renormalization fixed point are smoothly conjugated.}

\bigskip

Area-preserving maps at accumulation of period doubling are observed by several authors in the early 80's, see  \cite{DP,Hl,BCGG,Bo,CEK2,EKW1}. In \cite{EKW2} Eckmann, Koch and Wittwer introduced a period doubling renormalization scheme for area preserving maps and described the hyperbolic behavior of the renormalization operator in a neighborhood of a renormalization fixed point. In particular, they observed universality for maps at the accumulation of period doubling.

The maps in the universality class of the 
Eckmann-Koch-Wittwer renormalization fixed point are at the accumulation of period doubling. It was shown in  \cite{GJ1} that 
these maps do have a period doubling Cantor set and the Lyapounov exponents are at most zero. Moreover, for maps in a certain strong stable part of the renormalization fixed point, a space with finite codimension,
their period doubling Cantor sets are  at least bi-Lipschitz conjugated. In this paper we improve the conclusions of \cite{GJ1}: rigidity holds in the whole stable manifold of the Eckmann-Koch-Wittwer renormalization fixed point. Moreover, the conjugations are at least $C^{1+\alpha}$ with $\alpha=0.237$.
 
\bigskip

\noindent
{\bf Acknoledgements:} This work was started during a visit by the authors to the Institut Mittag-Leffler (Djursholm, Sweden) as part of the research program on ``Dynamics and PDEs''. The hospitality of the institute is gratefully acknowledged. The second author was funded by a postdoctoral fellowship from the Institut Mittag-Leffler. The authors would like to thank the referees for their very careful reading of the original manuscript.

\bigskip

%% file: prelim.tex
\section{Preliminaries}\label{prelim}

Given a domain $\mathcal{D}\subset \mathbb{C}^2$, let $D\subset \inter(\mathcal{D}\cap \mathbb{R}^2)$ be compactly contained in the real slice. Assume $(0,0)\in D$.
An {\it area-preserving map} $F:D\to F(D)\subset \mathbb{R}^2$ will mean a real symmetric map which has a holomorphic extension to $\mathcal{D}$ and is an exact symplectic diffeomorphism  onto its image with the following properties
\begin{itemize}
\item[1)] $F(0,0)=(0,0)$,
\item[2)]  $T \circ F \circ T=F^{-1}$,  where $T(x,u)=(x,-u)$ ({\it reversibility}),
\item[3)] $\partial_u x' \ne 0$ with $(x',u')=F(x,u)$ ({\it  twist condition}).
\end{itemize}
The collection of such maps is denoted by ${\rm Cons}(D)$. 

\bigskip

Our subject are maps at the accumulation of period doubling in a neighborhood of the Eckmann-Koch-Wittwer renormalization fixed point, \cite{EKW2}. We will need stronger and  additional  estimates than the spectral information obtained in \cite{EKW2}. These estimates were obtained in \cite{GJ2}. Although both \cite{EKW2} and \cite{GJ2} study a neighborhood of the same fixed point, the used constructions differ. We do need the estimates form \cite{GJ2} and summerize these results here. 

\bigskip

In \cite{GJ2} Gaidashev and Johnson construct simply connected domains   $D\subset \mathcal{D}\subset \mathbb{C}^2$ and adapt the renormalization scheme from \cite{EKW2}. This renormalization scheme is defined on a neighborhood $\text{Cons}_0(D)$ of $F_*\in 
\text{Cons}(D)$, where $F_*$ corresponds to the Eckman-Koch-Wittwer fixed point. Each map $F\in \text{Cons}_0(D)$ has a unique periodic orbit of period $2$. This orbit has one point $(p_F,0)$ with $p_F<0$.
The dependance
$$
F\mapsto p_F
$$
is analytic. Furthermore there are analytic function
$$
F\mapsto \lambda_F\in (-\infty,0)
$$
and
$$
F\mapsto \mu_F\in (0, \infty).
$$
The rescaling which will be used to renormalize is
$$
\Lambda_F:(x,u)\mapsto (\lambda_F x+p_F, \mu_F u).
$$
The renormalization operator $R:\text{Cons}_0(D)\to \text{Cons}(D)$ is defined by
$$
RF={\Lambda_F}^{-1}\circ F\circ F\circ \Lambda_F.
$$
The results from \cite{GJ2} which will be used in the sequel are collected in the following Theorem, Lemma, and Proposition. We will denote the $C^2$-norm of maps by $||. ||_{C^2}$.

\comm{
A map $F\in \text{Cons}(D)$ is called {\it renormalizable} if there exists a unique periodic point $(p_1,p_2)$ of period $2$ and a rectangle 
\begin{equation}\label{shape}
B_0=[p_1-a_1,p_1+a_1]\times [p_2-a_2,p_2+a_2]\subset [-1,1]^2
\end{equation}
such that
\begin{equation}\label{P1}
F^2(B_0)\cap B_0\ne \emptyset,
\end{equation}
and 
\begin{equation}\label{P2}
F(B_0)\cap B_0=\emptyset.
\end{equation}
The collection of renormalizable maps is denoted by $\text{Cons}_0(D)$. The renormalization operator $R:\text{Cons}_0(D)\to \text{Cons}(D)$ is defined as follows. Given a map $F\in \text{Cons}_0(D)$ let $B_0$ be the smallest rectangle with the properties (\ref{shape}), (\ref{P1}) and (\ref{P2}). Then
$$
RF={\Lambda_F}^{-1}\circ F\circ F\circ \Lambda_F,
$$
where $\Lambda_F: \mathbb{R}^2\to \mathbb{R}^2$ is the affine map
$$
\Lambda_F:(x,u)\mapsto (-\lambda_F x+h_1, \mu_F u+h_2)
$$
with $\Lambda_F(D)=B_0$ and $\lambda_F, \mu_F>0$. A map $F\in \text{Cons}_0(D)$ is infinitely renormalizable if $R^nF$ is defined for all $n\ge 1$.

The following Theorem is proved in \cite{GJ2}. See also Koch etc.
}

\begin{thm}\label{hyp} There exists $F_*\in \text{Cons}_0(D)$ such that
\begin{itemize}
\item[1)] $\text{Cons}_0(D)$ is a neighborhood of $F_*$.
\item[2)] $F_*$ is a hyperbolic fixed point of the renormalization operator.
\item[3)] $F_*$ has a one-dimensional local unstable manifold.
\item[4)] $F_*$ has a codimension one local stable manifold $W^s_{\text{loc}}(F_*)$.
\item[5)] $F^*$ has a codimension two local strong stable manifold 
$W^{ss}_{\text{loc}}(F_*)\subset W^s_{\text{loc}}(F_*)\subset \text{Cons}_0(D)$.
\item[6)] There exist a distance $d$ on $W^{ss}_{\text{loc}}(F_*)$, and $\nu<0.126$  such that for every $F, \tilde{F}\in W^{ss}_{\text{loc}}(F_*)$
$$
||F-\tilde{F}||_{C^2}\le d(F,\tilde{F})
$$
and 
$$
d(RF,R\tilde{F})\le \nu \cdot d(F,\tilde{F}).
$$
In particular,
$$
d(R^nF,F_*)\le   \nu^n\cdot d(F,F_*).
$$
\item[7)] The one dimensional family defined by $F_t=h_t^{-1}\circ F_*\circ h_t$, where $h_t:D\to h_t(D)\subset \mathbb{R}^2$ is the diffeomorphism defined by
$$
h_t(x,u)=(x+tx^2, \frac{u}{1+2tx}),
$$
with $|t|\le c_0$, is contained in the stable manifold $W^s_{\text{loc}}(F_*)$ and is transversal to the strong stable manifold $W^{ss}_{\text{loc}}(F_*)\subset W^s_{\text{loc}}(F_*)$ and intersects only in $F_*$.
\item[8)] The map $\Lambda_F$ depend analytically on $F$.
\end{itemize}
\end{thm}

For $F\in \text{Cons}_0(D)$ we will use the notation
$$
\psi_0=\psi_0(F)=\Lambda_{F}:D\to \mathbb{R}^2
$$
and
$$
\psi_1=\psi_1(F)=F\circ\Lambda_{F}:D\to \mathbb{R}^2
$$
In the next sections we will consider infinitely renormalizable maps. These are maps for which all $R^nF$, $n\ge 0$, are defined.  For such maps we will consider
$$
\Psi^n_w=\psi_{w_1}(F)\circ \psi_{w_2}(RF)\circ \cdots \circ \psi_{w_n}(R^{n-1}F):\Dom(R^nF)\to \Dom(F) ,
$$
with $w_k\in \{0,1\}$.

\begin{lem}\label{bob1} For $F\in \text{Cons}_0(D)$ there are analytically defined simply connected domains  $(p_F,0)\in B_0(F)\subset D$ and $B_1(F)\subset D$ such that
\begin{equation}\label{P2}
B_1(F)\cap B_0(F)=\emptyset,
\end{equation}
and 
\begin{equation}\label{P1}
F^2(B_0(F))\cap B_0(F)\ne \emptyset.
\end{equation}
Moreover,
$$
\psi_w(B_0(RF)\cup B_1(RF))\subset B_w(F),
$$
with $w\in \{0,1\}$.
\end{lem}

\begin{prop}\label{esti}  There exists $0<\theta_1<\theta_2<1$  such that for $F\in W^s_{\text{loc}}(F_*)$  we have
$$
\theta_1^4 \cdot |v|\le |D\Psi^4_w(x,u) v|\le \theta_2^4 \cdot  |v|
$$
for every $w\in \{0,1\}^4$ and $(x,u)\in B_0(F)\cup B_1(F)$. Moreover,
\begin{equation}\label{crucialestimate}
\frac{\theta_2\nu }{\theta_1}<1.
\end{equation}
\end{prop}

\begin{rem}
The following estimates are obtained in \cite{GJ2}.
\begin{equation}\label{theta1}
\theta_1\ge 0.061
\end{equation}
\begin{equation}\label{theta2}
\theta_2\le 0.249
\end{equation}
\begin{equation}\label{nu}
\nu< 0.126
\end{equation}
\begin{equation}\label{ratio}
\frac{\theta_2\nu }{\theta_1}< 0.515
\end{equation}
\end{rem}

\begin{rem}\label{optcond}
The estimate (\ref{crucialestimate}) plays a crucial role in the proof of the rigidity Theorem  \ref{rigthm}. The optimal way to describe this condition is in terms of $F_*$ and the derivative of renormalization $DR(F_*)$ at that point. Namely, let
$$
\sigma^*=\text{Spectral radius}(DR(F_*)|T^{ss}),
$$ 
where $T^{ss}$ is the tangent space to $W^{ss}_{\text{loc}}(F_*)$ at $F_*$. Morever, let
$$
\theta^*_1=\inf_{w\in \{0,1\}^\infty, (x,u)\in D, ||v||=1} \{
\liminf_{ n\to \infty}  ||D\Psi^n_{w^n}(F_*)(x,u)v||^{\frac{1}{n}}
\},
$$
where $w^n$ is the truncation of the infinite word $w$ to its first $n$ letters. And similarly,
let
$$
\theta^*_2=\sup_{w\in \{0,1\}^\infty, (x,u)\in D, ||v||=1} \{
\limsup_{ n\to \infty}  ||D\Psi^n_{w^n}(F_*)(x,u)v||^{\frac{1}{n}}
\}.
$$
The numbers $\theta^*_1$ and $\theta^*_2$ reflect geometrical properties of the period doubling Cantor set of $F_*$, see \S\ref{cantor}. The rigidity Theorem can be proved under the condition
$$
\frac{\theta^*_2 \sigma^*}{\theta^*_1}<1.
$$
The derivative of renormalization at the fixed point is a compact operator. In particular,  rigidity can be proved on a finite codimension subspace where the contraction is strong enough. The numerical estimates from \cite {GJ2} show that only the weakest stable direction is not strong enough. Luckily, this weakest stable direction corresponds to a one-dimensional family of analytically conjugated maps. The authors do not have a conceptual explanation for this coincidence.
\end{rem}

The Appendix is used to present some elementary estimates on the error term of affine approximations of diffeomorphisms. In this section we used coordinates $(x,u)$ to denote a point. In the sequel we will not need anymore the coordinates and denote points simply by $x$, etc.

%% file: cantorset.tex
\section{The Cantor Set}\label{cantor}

\begin{figure}[htbp]
\begin{center}
\psfrag{F}[c][c] [0.7] [0] {\Large $F$}
\psfrag{F2}[c][c] [0.7] [0] {\Large $RF$}
\psfrag{F3}[c][c] [0.7] [0] {\Large $RF^2$}
\psfrag{Fn}[c][c] [0.7] [0] {\Large $R^nF$}
\psfrag{B0}[c][c] [0.7] [0] { $B_0$}
\psfrag{B1}[c][c] [0.7] [0] { $B_1$}
\psfrag{p01}[c][c] [0.7] [0] {\Large $\psi_0^1$} 
\psfrag{p11}[c][c] [0.7] [0] {\Large $\psi_1^1$} 
\psfrag{p02}[c][c] [0.7] [0] {\Large $\psi_0^2$} 
\psfrag{p12}[c][c] [0.7] [0] {\Large $\psi_1^2$} 
\psfrag{p03}[c][c] [0.7] [0] {\Large $\psi_0^3$} 
\psfrag{p13}[c][c] [0.7] [0] {\Large $\psi_1^3$} 
\pichere{1.0}{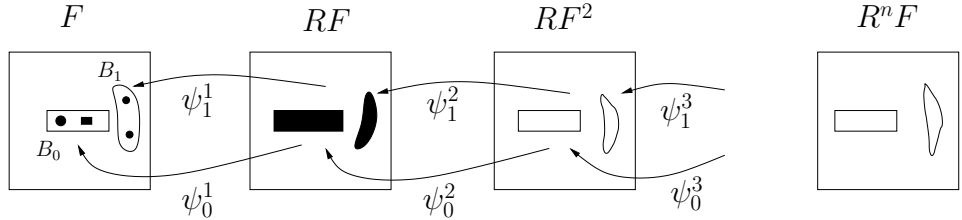} 
\caption{The Renormalization Microscope} \label{micro}
\end{center}
\end{figure}

In this section we construct the invariant Cantor set for infinitely renormalizable maps. This cantor set is intrinsically related to the renormalization process. As in dimension one, it is a Cantor set on which the map acts like the dyadic adding machine. The construction is similar to the construction of the limit set of an iterated function system. As a matter of fact, the Cantor set of the renormalization fixed point $F_*$, is the limit set of an iterated function system, see Figure \ref{micro}.

The stable manifold $W^s(F_*)\subset \text{Cons}_0(D)$, also called {\it the universality class},  of $F_*$, consists of the maps with $R^nF\to F_*$.
 Trough out this section we will consider a fixed map $F\in W^s(F_*)$.

We will use the notation
$$
\psi_0^n=\Lambda_{R^{n-1}F}:D\to \mathbb{R}^2
$$
and
$$
\psi_1^n=R^{n-1}F\circ\Lambda_{R^{n-1}F}:D\to \mathbb{R}^2
$$
Observe,
$$
R^nF=(\psi_0^n)^{-1}\circ R^{n-1}F\circ \psi_0^n.
$$
The convergence of $R^nF\to F_*$ and Theorem \ref{hyp}(8) immediately implies

\begin{lem}\label{limitpsi} For every $F\in W^s(F_*)$
$$
\lim_{n\to\infty} ||\psi_{0,1}^n(F)-\psi_{0,1}(F_*)||_{C^2}=0.
$$
\end{lem}

\begin{rem} \label{degeneration}Lemma \ref{limitpsi} shows a crucial difference between conservative and dissipative  infinitely renormalizable maps. In the conservative case, the rescaling maps converge to (non-degenerated) diffeomorphisms. In the dissipative case, the corresponding rescaling $\psi_1^n$ converge to a degenerated map, a map with one-dimensional image. This degeneration is at the heart of the difficulties of the non-rigidity phenomena observed in dissipative maps, see \cite{CLM}. The degeneration occurs in a universal manner. This universal degeneration is responsible for the probabilistic nature of geometry in dimension two where one observes Probabilistic Universality and Probabilistic Rigidity, see \cite{LM3}.
\end{rem}

\begin{rem}\label{smoothness} The convergence in Lemma \ref{limitpsi}
holds in any $C^k$-distance, $k\ge 0$. We only need control of the lower order derivatives.
\end{rem}

\begin{rem}\label{constr} The construction of the Cantor set in the conservative case is exactly the same as in the dissipative case. The difference lies in the asymptotic behavior of the rescalings, see  remark \ref{degeneration}.
\end{rem}

Let
$$
\Psi^2_{00}= \psi^1_0\circ \psi^2_0, \quad \Psi^2_{01}= \psi^1_0\circ
\psi^2_1, \quad \Psi^2_{10}=\psi^1_1\circ\psi^2_0,\quad \dots.
$$
and, proceeding this way, construct, for any  $w=(w_1, \dots, w_n)\in\{0,1\}^n$, $n\ge 1$,
the maps
$$
 \Psi^n_w = \psi^1_{w_1}\circ\dots\circ \psi^n_{w_n}:\Dom(R^nF)\to \Dom(F).
 $$
The notation $\Psi^n_w(F)$ will also be used to emphasize the map
under consideration. The following Lemma follows directly from Proposition \ref{esti}.

\begin{lem}\label{contracting}
For every $F\in W^s(F_*)$ there exists $C>0$ such that 
for any word $w\in \{0,1\}^n$, $n\ge 1$,
$$
\| D\Psi^n_w\|\leq  C \theta_2^n
$$
where $\theta_2<1$ is given in Proposition \ref{esti} and (\ref{theta2}).
\end{lem}

Define the {\em pieces} of the  $n^{\text{th}}$-{\em level} or $n^{\text{th}}$-{\em scale} as follows. They are closed
topological disks.  For each $w\in \{0,1\}^n$ let
$$
B^n_{w0}\equiv B^n_{w0}(F) = \Psi^{n}_w (B_{0}(R^nF))
$$
and
$$
B^n_{w1}\equiv B^n_{w1}(F) = \Psi^{n}_w (B_{1}(R^nF)).
$$

The set of words $\{0,1\}^n$ can be viewed as the additive group of residues $\mod 2^n$ by letting 
$$ w
\mapsto\sum_{k=0}^{n-1} w_{k+1} 2^k.
 $$
Let $p\colon \{0,1\}^n\ra \{0,1\}^n$ be
the operation of adding 1 in this group.

\begin{lem}\label{permute}
For every $F\in W^s(F_*)$ and $n\ge 1$
\begin{enumerate}
\item[1)] The above families of pieces are nested:
$$
   B^n_{w\nu}\subset B^{n-1}_w, \quad w\in \{0,1\}^{n}, \ \nu\in \{0,1\}.
$$
\item[2)] The pieces $B^n_w$, $w\in \{0,1\}^{n+1}\setminus \{1^{n+1}\}$, are pairwise
disjoint.
\item[3)] Under $F$, the pieces are permuted as follows.
$$
F(B^n_w) = B^n_{p(w)},
$$ 
unless $p(w) = 0^{n+1}$. If
$p(w)=0^{n+1}$, then $F(B^n_w) \cap B^n_{0^{n+1}}\ne \emptyset$.
\end{enumerate}
\end{lem}

\begin{proof} The first assertion holds by construction: Let $w\in \{0,1\}^{n-1}$, $w_n,v\in \{0,1\}$, and use Lemma \ref{bob1},
$$
\begin{aligned}
B^n_{ww_nv}&= \Psi^{n}_{ww_n } (B_v(R^nF))= \Psi^{n-1}_{w}\circ
\psi^{n}_{ w_n} (B_v(R^nF))\\
&\subset \Psi^{n-1}_{w}(B_{w_n}(R^{n-1}F))
=B^{n-1}_{ww_n}(F).
\end{aligned}
$$
The second follows by induction.
 For all maps under consideration we have the disjointness property (\ref{P2}) which says, 
 that $B^1_0(F)$ and $B^1_1(F)$ are disjoint. Assume that the pieces of the
 $n^{th} $ generation
are disjoint for all maps under consideration. This implies that
the pieces $B^{n+1}_{0w}\subset B^1_0$, $w\in \{0,1\}^{n+1}$,  are pairwise disjoint,
as they are images of the disjoint pieces $B^{n}_{w }(RF)$  by the map $\psi_1^0$. 
Applying $F$, we see that the pieces $B^{n+1}_{1w}\subset B^1_1$, $w\in W^n$,
 are pairwise disjoint as well.
The assertion follows because  $B_1^1$ and $B_0^1$ are also disjoint.

\ssk Let us inductively check the third assertion. For $n=1$, we
have:
$$
B^1_1= F (B^1_0) \text{ and } F(B^1_1) \cap B^1_0\ne \emptyset.
$$
Consider now  the pieces $B^n_w(RF)$, $w\in \{0,1\}^{n+1}$, of level $n$
for $RF$. Assume inductively that they are permuted by $RF$ as
required. Then the pieces $B^{n+1}_{0 w} = \psi^1_0 (B^n_w(RF))$
 are permuted in the same fashion under $F^2$.
Moreover, $B^{n+1}_{1w} = \psi^1_1 (B^n_w(RF)) = F(B^{n+1}_{0w})$, 
and the conclusion follows.
\end{proof}

Furthermore, Lemma \ref{contracting} implies:

\begin{lem}\label{boxes shrink}
For every $F\in W^s(F_*)$ there exists $C>0$ such that
for all $w\in \{0,1\}^{n+1}$, $ \diam B^n_w\leq C \theta_2^n$.
\end{lem}

Let
$$
\OO\equiv \OO_F= \bigcap_{n=1}^\infty \bigcup_{w\in \{0,1\}^n} B^n_w.
$$
Let us also consider the {\it diadic group} $\displaystyle { \{0,1\}^\infty = \lim_{\leftarrow} \{0,1\}^n}$.
The elements of $\{0,1\}^\infty$ are infinite sequences $(w_1 w_2\dots)$ of symbols $0$ and $1$
that can be also represented as formal power series
$$
   w \mapsto\sum_{k=0}^\infty w_{k+1} 2^k.
$$
The integers $\Z$ are embedded into $\{0,1\}^\infty$ as finite series.
The {\it adding machine} $p: \{0,1\}^\infty\ra \{0,1\}^\infty$ 
is the operation of adding $1$ in this group.   
The discussion above yields:

\begin{thm}\label{adding machine} For every $F\in W^s(F_*)$
the set $\OO_F$ is an invariant Cantor set. The map $F$ acts on 
$\OO_F$ as the dyadic adding machine $p: \{0,1\}^\infty\ra \{0,1\}^\infty$.
The conjugacy between $p$ and $F|\OO_F$ is given by the following homeomorphism $h_F: \{0,1\}^\infty\ra \OO_F$: 
\begin{equation}\label{conjh}
     h_F: w = (w_1w_2\dots) \mapsto \bigcap_{n=1}^\infty B^n_{w_1\dots w_{n+1}}. 
\end{equation}
 Furthermore, $\OO_F$ has Lebesgue measure zero with
$$
\text{HD}(\OO_F)\le -\frac{\log 2}{\log \theta_2}\le 0.795. 
$$
\end{thm}

The invariant Cantor sets $\OO_F$ are the counterpart of the period doubling Cantor sets in one-dimensional dynamics and strongly dissipative H\'enon-like maps, see \cite{CLM, GST, Mi}. The dynamics of $F$ restricted to $\OO_F$ is conjugated to the adding machine. The adding machine is uniquely ergodic. Let $\mu$ be the unique invariant measure of $F$ restricted to $\OO_F$:
$$
\mu(B^n_{w})=\frac{1}{2^{n+1}}.
$$

\begin{thm}\label{expo} The measure $\mu_F$ of every map $F\in W^s(F_*)$ has a single characteristic exponent, $\chi=0$.
\end{thm}

\begin{proof} The largest characteristic exponent is denoted by $\chi$. Let $F_n$ be the $n$-th renormalization of $F$. This map restricted to $B_0(R^nF)\cup B_1(R^nF)$ is smoothly
conjugate to the restriction of $F^{2^n}$ to the piece
$B^n_{0^n0}\cup B^n_{0^n1}$. Let $\mu_n$ be the normalized restriction of $\mu$ to
$B^n_{0^n0}\cup B^n_{0^n1}$, and let $\nu_n$ be the invariant measure on the Cantor set of
$F_n$. Note that these two measures are preserved by the
conjugacy. Then 
$$ 2^n \chi = \chi (F^{2^n}| B^n_{0^n0}\cup B^n_{0^n1},\, \mu_n) =
\chi(F_n , \nu_n) \leq
\int \log \|DF_n\|\,  d\nu_n \leq C,
$$
since the maps $F_n\to F_*$ have uniformly bounded $C^1$-norms.

Hence, $\chi\leq 0$. If $\chi<0$, both characteristic exponents of
$F$ would be negative. This contradicts  the
relation $\chi+ \chi_- = 0$ which holds because the map preserves the Lebesgue measure.
\end{proof}

%% file: rigidity.tex
\section{Rigidity}\label{rigidity}

Let $ \OO\subset \mathbb{R}^2$. A map $h:\OO\to \mathbb{R}^2$
is {\it differentiable} at $x_0\in \OO$ if
there exists a linear map $Dh(x_0): \mathbb{R}^2\to \mathbb{R}^2$ such that for 
 $x\in \OO$
$$
h(x)=h(x_0)+Dh(x_0)(x-x_0)+o(|x-x_0|).
$$
If the map $x\mapsto Dh(x)$ is $C^{\alpha}$ then we say that  $h$ is a $C^{1+\alpha}$-map and if it  is Lipschitz ($\alpha=1$) we say that $h$ is $C^{1+\text{Lip}}$. Observe, the composition of two $C^{1+\alpha}$-maps is again $C^{1+\alpha}$.

A bijection $h:\OO\to h(\OO)\subset \mathbb{R}^2$ is a {\it $C^{1+\alpha}$-diffeomorphism} if there exists a map $h^{-1}: h(\OO)\to \mathbb{R}^2$
with $h^{-1}\circ h=\text{id}$ and both maps are $C^{1+\alpha}$.

There are many conjugations between $\OO_F$ and $\OO_{\tilde{F}}$
with $F, \tilde{F}\in W^s(F_*)$. However, we will only consider conjugations $h:\OO_F\to \OO_{\tilde{F}}$ defined
by
$$
h=h_{\tilde{F}}\circ h_F^{-1},
$$
where $h_F$ and $h_{\tilde{F}}$ are the conjugations from $\{0,1\}^\infty$ to $\OO_F$ and $\OO_{\tilde{F}}$ respectively, see Theorem \ref{adding machine}. 

\begin{defn}\label{rigid} The Cantor set $\OO_F$ is $rigid$ if  for some $\alpha>0$, the conjugation $h:\OO_F\to \OO_{F_*}$ is a $C^{1+\alpha}$-diffeomorphsim. Notation,
$
\OO_F = \OO_{F_*} \mod (C^{1+\alpha}).
$
\end{defn}

Let
\begin{equation}\label{alpha0}
\alpha_0=\frac{\ln \theta_2 \nu}{\ln \theta_1}-1>0.237.
\end{equation}

\begin{thm}\label{rigthm}  The Cantor set $\OO_F$, with $F\in W^{s}_{\text{loc}}(F_*)$, is rigid. Namely, 
$$
\OO_F =\OO_{F_*} \mod(C^{1+\alpha})
$$
for every $0<\alpha<\alpha_0$.
\end{thm}

The proof of the Rigidity Theorem \ref{rigthm}  consists of two parts. The main part, Proposition \ref{core} treats the rigidity question for maps in the local strong stable manifold. The other reduces the rigidity question to this Proposition, using that the local weak stable manifold is part of a smooth conjugacy class, see Lemma \ref{famtrans}.

\begin{lem}\label{RnFconj} If $
\OO_{R^nF} =\OO_{F_*} \mod(C^{1+\alpha})
$
then
$
\OO_F =\OO_{F_*} \mod(C^{1+\alpha}).
$
\end{lem}

\begin{proof} Assume that the conjugation $\hat{h}_{n}: \OO_{R^nF}\to \OO_{F_*}$ is a $C^{1+\alpha}$-diffeomorphism. Use the notation
$\OO^k_w=B^k_w\cap \OO_F$ (To indicate the map under consideration
we will also use
$\OO^k_w(F)=B^k_w\cap \OO_F$) . Observe, for every word $w\in \{0,1\}^{k+1}$ one has 
\begin{equation}\label{hnpieces}
\hat{h}_n(\OO^k_w(R^nF))=\OO^k_w(F_*).
\end{equation}
Define the following $C^{1+\alpha}$-diffeomorphism $h_n: \OO(F)\to \OO(\tilde{F})$. For $w\in \{0,1\}^n$, $v\in \{0,1\}$ let
\begin{equation}\label{hn}
h_n|\OO_{wv}^n=\Psi^n_w(F_*) \circ \hat{h}_{n}\circ (\Psi^n_w(F))^{-1}: \OO_{wv}^n(F)\to \OO_{wv}^n(\tilde{F}).
\end{equation}
From (\ref{hnpieces}) we get a similar property for $h_n$: for every word $w\in \{0,1\}^{k+1}$, $k\ge n$,
$$
h_n(\OO^k_w(F))=\OO^k_w(F_*).
$$
Let $x\in \OO_F$, say
$
x=\bigcap_k \OO^k_{w^k}(F),
$
with $w^k\in \{0,1\}^k$.
Observe, 
$$
F(x)=\bigcap_k \OO^k_{p(w^k)}(F)
$$
and
$$
h(x)=\bigcap_k \OO^k_{w^k}(F_*).
$$
Hence,
$$
F_*(h(x))=\bigcap_k \OO^k_{p(w^k)}(F_*)
$$
and
$$
h(F(x))=\bigcap_k \OO^k_{p(w^k)}(F_*).
$$
Indeed, the diffeomorphism is a conjugation, $h\circ F=F_*\circ h$.
\end{proof}

Consider the family of diffeomorpsims  $h_t:D\to h_t(D)\subset \mathbb{R}^2$ given in Theorem \ref{hyp}(7). 

\begin{lem}\label{famtrans} There exists $c>0$ such that for any  $F\in W^{s}_{\text{loc}}(F_*) $ close enough to $F_*$ the family $F_t=h_t^{-1}\circ F\circ h_t$, $|t|\le c$, has a unique intersection with  the strong stable manifold $W^{ss}_{\text{loc}}(F_*)$.
\end{lem}

\begin{proof} This Lemma follows immediately from the transversality property mentioned in Theorem \ref{hyp}(7).
\end{proof}

\begin{lem}\label{friend} For every $F\in W^{s}(F_*)$ there exists  
$\tilde{F}\in W^{ss}_{\text{loc}}(F_*)$  such that for some $n\ge 1$,  the conjugation $h:\OO_{R^nF}\to \OO_{\tilde{F}}$ is a $C^{1+\text{Lip}}$-diffeomorphism.
\end{lem}

\begin{proof} For $n\ge 1$ large enough, $R^nF\in W^s_{\text{loc}}(F_*)$ is close enough to $F_*$ to apply Lemma \ref{famtrans}.
\end{proof}

\begin{rem} Observe, the conjugation in lemma \ref{friend} is much better than $C^{1+\text{Lip}}$. It is in fact the restriction of a rational map.
\end{rem}

\begin{prop}\label{core}  The Cantor set $\OO_F$, with $F\in W^{ss}_{\text{loc}}(F_*)$, is rigid. Namely, 
$$
\OO_F =\OO_{F_*} \mod(C^{1+\alpha})
$$
for every $0<\alpha<\alpha_0$.
\end{prop}

The proof of this key Proposition will be presented in a series of Lemmas.
Fix two maps $F,\tilde{F}\in W^{ss}_{\text{loc}}(F_*)$.  The constants which will appear in the following analysis are independent of these maps. They only depend on the size of $W^{ss}_{\text{loc}}(F_*)$. 

\bigskip

Consider the conjugation $h:\OO_F\to \OO_{\tilde{F}}$ defined
by
$$
h=h_{\tilde{F}}\circ h_F^{-1},
$$
where $h_F$ and $h_{\tilde{F}}$ are the conjugations from the adding machine $\{0,1\}^\infty$ to $\OO_F$ and $\OO_{\tilde{F}}$ respectively, see Theorem \ref{adding machine}. We will show that this map is $C^{1+\alpha}$, for any $\alpha>0$
satisfying
\begin{equation}\label{alpha}
\frac{\theta_2  \nu }{\theta_1^{1+\alpha}}< 1
\end{equation}
or equivalently $0<\alpha<\alpha_0$, see (\ref{alpha0}).

\bigskip

Proposition \ref{esti} has an unfortunate form. Ideally, we would have 
$$
\theta_1 \cdot |v|\le |D\psi_{0,1}(x,u) v|\le \theta_2 \cdot  |v|.
$$
However, this only holds on a small neighborhood of the Cantor sets. This inconvenience forces us to work in "steps of four".

\bigskip 

\begin{figure}[htbp]
\begin{center}
\psfrag{F}[c][c] [0.7] [0] {\Large $F$}
\psfrag{Tf}[c][c] [0.7] [0] {\Large $\tilde{F}$}
\psfrag{Fn}[c][c] [0.7] [0] {\Large $R^{4n}F$}
\psfrag{TFn}[c][c] [0.7] [0] {\Large $R^{4n}\tilde{F}$}
\psfrag{Bnw}[c][c] [0.7] [0] { $B^{4n}_{w0}$}
\psfrag{TBnw}[c][c] [0.7] [0] { $\tilde{B}^{4n}_{w0}$}
\psfrag{hn}[c][c] [0.7] [0] {\Large $h_n$} 
\psfrag{id}[c][c] [0.7] [0] {\Large $\id$} 
\psfrag{psnw}[c][c] [0.7] [0] {\Large $\Psi^{4n}_w$} 
\psfrag{tpswn}[c][c] [0.7] [0] {\Large $\tilde{\Psi}^{4n}_w$} 
\psfrag{O}[c][c] [0.7] [0] {\Large $\OO_F$}
\psfrag{TO}[c][c] [0.7] [0] {\Large $\OO_{\tilde{F}}$}
\psfrag{b0rf}[c][c] [0.7] [0] { $B_0$}  
\psfrag{b0rtf}[c][c] [0.7] [0] {$\tilde{B}_0$}
\pichere{0.9}{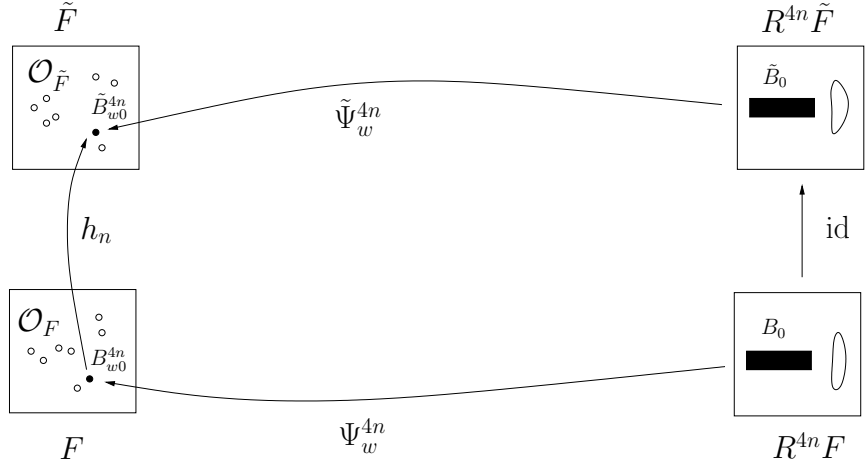} 
\caption{The Aproximate Conjugation $h_n$} \label{aprox}
\end{center}
\end{figure}

Let $\OO_n=\bigcup_{w\in \{0,1\}^{4n+1}} B^{4n}_w(F)$ and similarly define $\tilde{\OO}_n$ for $\tilde{F}$   (In the sequel we will use  a tilde to indicate whether an object refers to $F$ or $\tilde{F}$). 
The map $h_n:\OO_n\to \tilde{\OO}_n$ is defined  by
$$
h_n|B^{4n}_{wv}=\tilde{\Psi}^{4n}_w\circ (\Psi^{4n}_w)^{-1},
$$
where $w\in \{0,1\}^{4n}$ and $v\in \{0,1\}$, $\Psi^{4n}_w=\Psi^{4n}_w(F)$ and $\tilde{\Psi}^{4n}_w=\Psi^{4n}_w(\tilde{F})$. This map $h_n$  is an approximate conjugation, see Figure \ref{aprox}. Compare this map with the conjugation (\ref{hn}). The conjugation $\hat{h}_{n}$ between the $n^{th}$-renormalization is replaced here by the identity. To indicate the maps under consideration we will also use the notation $h_n(F,\tilde{F})$.

Observe, $h_n|\OO_F\to h$, where $h:\OO_F\to \OO_{\tilde{F}}$ is the conjugation. The following Lemmas will describe this convergence in more detail. The proofs are by induction. The idea of this induction is the following. 
The maps $h_n$ are quasi conjugations up to the $4n^{th}$-scale. We use the uniform control of these maps at the $4n^{th}$-scale for the renormalizations and push it down to the next scale of the original systems. The control of the derivative of the quasi-conjugations, see Lemma \ref{c1conv}, Lemma \ref{c2conv},  and Lemma \ref{caDh} relies in an essential way on the fact that the renormalizations $R^4F$ and $R^4\tilde{F}$ of $F$ and $\tilde{F}$ are strictly closer, see Theorem \ref{hyp}(6), than the original systems. Namely, this implies that their quasi-conjugation at the $4n^{th}$-scale is strictly better than the quasi-conjugation at the $n^{th}$-scale of the original systems. This effect causes the contraction in for example the estimate (\ref{contr}).

\begin{lem}\label{c0conv} 	There exists a constant $C>0$ such that 
$$
|h_{n+1}(x)-h_n(x)|\le C\cdot (\theta_2\nu)^{4n}\cdot d(F, \tilde{F}),
$$
with $x\in \OO_{n+1}$. In particular,
$$
|h_n(x)-x|\le C\cdot d(F, \tilde{F})
$$
\end{lem} 

\begin{proof} The proof will be by induction. Observe, $h_0(x)=x$ and
$$
h_1(x)=\tilde{\Psi}^4_w \circ (\Psi^4_w)^{-1}(x),
$$
when $x\in B^4_w$ with $w\in \{0,1\}^{4+1}$. Because the maps $\Psi^4_w(F)$ depend analytically on $F$, see Theorem \ref{hyp}(8), we get
$$
|h_1(x)-h_0(x)|\le C\cdot d(F, \tilde{F}).
$$
Assume the Lemma holds for $ n$. Let $\hat{h}_n=h_n(R^4F,R^4\tilde{F})$. Given a point in the domain of $h_{n+2}$, say $x\in B^{4(n+2)}_{wv}$ with $w\in \{0,1\}^4$ and $v\in \{0,1\}^{4(n+1)+1}$. 
Then
$$
h_{k+1}(x)=\tilde{\Psi}^4_w \circ \hat{h}_k \circ (\Psi^4_w)^{-1}(x),
$$
with either $k=n$ or $k=n+1$. Also observe that
$(\Psi^4_w)^{-1}(B^{4(n+2)}_{wv})\subset B^{4(n+1)}_{v}\subset B_0(R^4F)\cup B_1(R^4F)$.
This will allow us to apply Proposition \ref{esti}.
Let $x_-= (\Psi^4_w)^{-1}(x)$. Use Proposition \ref{esti}, Theorem \ref{hyp}(6), and the induction hypothesis, to obtain
$$
\begin{aligned}
|h_{n+2}(x)-h_{n+1}(x)|&\le
| \tilde{\Psi}^4_w \circ \hat{h}_{n+1}(x_-) -\tilde{\Psi}^4_w \circ \hat{h}_n(x_-)|\\
&\le \theta_2^4\cdot | \hat{h}_{n+1}(x_-) - \hat{h}_n(x_-)|\\
&\le \theta_2^4\cdot C\cdot (\theta_2\nu)^{4n}\cdot d(R^4F, R^4\tilde{F})\\
&\le C \cdot (\theta_2\nu)^{4(n+1)}\cdot d(F, \tilde{F}).
\end{aligned}
$$
The proximity to identity follows from $h_0=\id$.
\end{proof}

\comm{
\begin{proof} Let $x\in B^{4(n+1)}_{wv}$, with $w\in \{0,1\}^{4n}$, $v\in \{0,1\}^4$. The points $x_n,x_{n+1}\in D$ are such that $x=\Psi^{4n}_w(x_n)$ and $x_n=\Psi^4_{v}(R^{4n}F)(x_{n+1})$ .  Observe,
$$
h_n(x)=\tilde{\Psi}^{4n}_w(x_n).
$$
Similarily, let $\tilde{z}_n=\tilde{\psi}^{4}_v(R^{4n}\tilde{F})(x_{n+1})$
and observe,
$$
h_{n+1}(x)=\tilde{\Psi}^{4n}_w(\tilde{z}_n).
$$
Moreover,
$$
|\tilde{z}_n-\tilde{x}_n|\le C\cdot \nu^{4n}\cdot d_{C^2}(F,\tilde{F}).
$$
This estimate  follows by applying Theorem \ref{hyp}(8) 
$$
d_{C^0}(\Psi^{4}_v(R^{4n}\tilde{F}), \Psi^4_{v}(R^{4n}F))\le Cd_{C^2}(R^{4n}\tilde{F}, R^{4n}F)\le C\nu^{4n}\cdot d_{C^2}(\tilde{F}, F).
$$
Now apply Proposition \ref{esti}
$$
\begin{aligned}
|h_{n+1}(x)-h_n(x)|&\le \max_z |D\tilde{\Psi}^{4n}_w(z)|\cdot |\tilde{z}_n-\tilde{x}_n|\\
&\le C\cdot (\theta_2\nu)^{4n}\cdot d_{C^2}(F,\tilde{F}).
\end{aligned}
$$
The proximity to identity follows from $h_0=\id$.
\end{proof}
}

Observe that the maps $h_n$ are differentiable, in fact analytic.

\begin{figure}[htbp]
\begin{center}
\psfrag{F}[c][c] [0.7] [0] {\Large $F$}
\psfrag{Tf}[c][c] [0.7] [0] {\Large $\tilde{F}$}
\psfrag{Fn}[c][c] [0.7] [0] {\Large $R^{4}F$}
\psfrag{TFn}[c][c] [0.7] [0] {\Large $R^{4}\tilde{F}$}
\psfrag{Bnw}[c][c] [0.7] [0] { $B^{4}_{w1}$}
\psfrag{TBnw}[c][c] [0.7] [0] { $\tilde{B}^{4}_{w1}$}
\psfrag{hn}[c][c] [0.7] [0] {\Large $h_{n+1}$} 
\psfrag{nhat}[c][c] [0.7] [0] {\Large $\hat{h}_{n}$} 
\psfrag{psnw}[c][c] [0.7] [0] {\Large $(\Psi^{4}_w)^{-1}$} 
\psfrag{tpsnw}[c][c] [0.7] [0] {\Large $\tilde{\Psi}^{4}_w$} 
\psfrag{b1}[c][c] [0.7] [0] {\Large $B_1(R^4F)$}
\psfrag{tb1}[c][c] [0.7] [0] {\Large $B_1(R^4\tilde{F})$}  
\pichere{0.6}{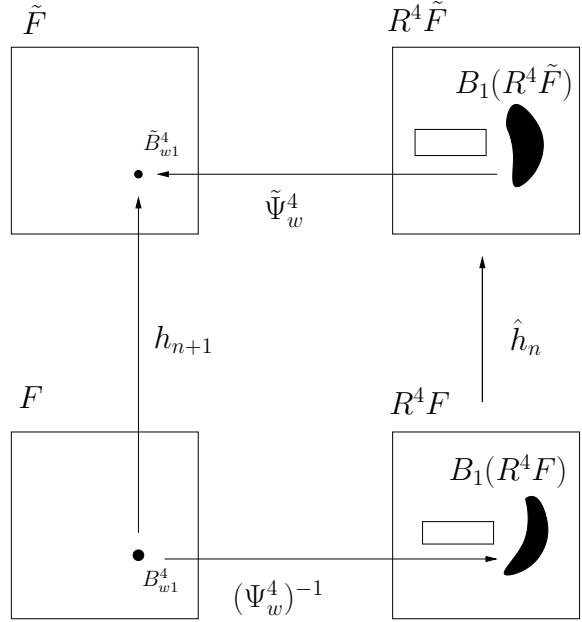} 
\caption{The Induction Step for  $h_{n+1}$} \label{induc}
\end{center}
\end{figure}

\begin{lem}\label{c1conv} 	There exists a constant $C>0$ such that 
$$
|Dh_{n+1}(x)-Dh_n(x)|\le C\cdot \left(\frac{\theta_2\nu}{\theta_1} \right)^{4n}\cdot d(F, \tilde{F}),
$$
with $x\in \OO_{n+1}$. In particular,
$$
|Dh_n(x)-\id|\le C\cdot d(F, \tilde{F}).
$$
\end{lem}

\begin{proof} The proof will be by induction. Observe, $Dh_0(x)=\id$ and
$$
h_1(x)=\tilde{\Psi}^4_w \circ (\Psi^4_w)^{-1}(x),
$$
when $x\in B^4_{wv}$ with $w\in \{0,1\}^{4}$ and $v\in \{0,1\}$. Because the maps $\Psi^4_w(F)$ depend analytically on $F$, see Theorem \ref{hyp}(8), we get
$$
|Dh_1(x)-Dh_0(x)|\le C\cdot d(F, \tilde{F}).
$$
Assume the Lemma holds for $0,1,2,\dots,n$ with constant $C_n>0$. This implies
\begin{equation}\label{Dhbound}
|Dh_n(x)-\id |\le C\cdot C_n \cdot d(F, \tilde{F}),
\end{equation}
where we used $Dh_0=\id$.

Let $\hat{h}_n=h_n(R^4F,R^4\tilde{F})$. Given a point in the domain of $h_{n+2}$, say $x\in B^{4(n+2)}_{wv}$ with $w\in \{0,1\}^4$ and $v\in \{0,1\}^{4(n+1)+1}$. 
Then
$$
Dh_{k+1}(x)=D\tilde{\Psi}^4_w(\hat{x}_k) \circ D\hat{h}_k(x_-) \circ D(\Psi^4_w)^{-1}(x),
$$
where $x_-= (\Psi^4_w)^{-1}(x)$ 
and $\hat{x}_k=\hat{h}_k(x_-)$ and either $k=n$ or $k=n+1$. Also observe that
$(\Psi^4_w)^{-1}(B^{4(n+2)}_{wv})\subset B^{4(n+1)}_{v}\subset B_0(R^4F)\cup B_1(R^4F)$.
This will allow us to apply Proposition \ref{esti}.
Let $d=d(F, \tilde{F})$.
Use Proposition \ref{esti}, Theorem \ref{hyp}(6), the induction hypothesis,
Lemma \ref{c0conv}, and (\ref{Dhbound}) to obtain
$$
\begin{aligned}
&|Dh_{n+2}(x)-Dh_{n+1}(x)| \\
\le &| D\tilde{\Psi}^4_w(\hat{x}_{n+1})  D\hat{h}_{n+1}(x_-)  D(\Psi^4_w)^{-1}(x)-
D\tilde{\Psi}^4_w(\hat{x}_n)  D\hat{h}_n(x_-) D(\Psi^4_w)^{-1}(x)| \\
\le &|D\tilde{\Psi}^4_w(\hat{x}_n)|
\cdot | D\hat{h}_{n+1}(x_-) - D\hat{h}_n(x_-)|\cdot |D(\Psi^4_w)^{-1}(x)|+\\
&|D\tilde{\Psi}^4_w(\hat{x}_{n+1})-D\tilde{\Psi}^4_w(\hat{x}_n) |\cdot
|D\hat{h}_{n+1}(x_-) |\cdot |D(\Psi^4_w)^{-1}(x)| 
\\
\le & \left(\frac{\theta_2}{\theta_1}\right)^4\cdot C_n\cdot \left(\frac{\theta_2\nu}{\theta_1}\right)^{4n}\cdot d(R^4F, R^4\tilde{F})+
C\cdot |\hat{x}_{n+1}-\hat{x}_{n}|\cdot |D\hat{h}_{n+1}(x_-) |\\
\le &  \text{   } C_n\cdot \left(\frac{\theta_2\nu}{\theta_1}\right)^{4(n+1)}\cdot d+
C\cdot (\theta_2\nu)^{4n} \cdot (1+C\cdot C_n\cdot d)\cdot d\\
\end{aligned}
$$
This implies
$$
C_{n+1}\le C_n+
C\cdot \theta_1^{4n} \cdot (1+C\cdot C_n).
$$
In particular, the sequence $C_n$ is bounded. The Lemma is proved.
\end{proof}

\begin{lem}\label{c2conv} 	For every $0<\alpha<\alpha_0$ there exists a constant $C>0$ such that 
$$
|h_{n}(y)-h_n(x)-Dh_n(x)(y-x)|\le C  \cdot d(F, \tilde{F})\cdot |y-x|^{1+\alpha},
$$
where $x,y\in \OO_{n}$.
\end{lem}

\begin{proof} Fix $<\alpha<\alpha_0$. 
Consider the collection of diffeomorphisms 
$$
\{\Psi^4_w(F)| F\in  W^{ss}_{\text{loc}}(F_*) \text{ and } w\in \{0,1\}^4\}
$$
and
$$
\{ (\Psi^4_w(F))^{-1}| F\in  W^{ss}_{\text{loc}}(F_*) \text{ and } w\in \{0,1\}^4\}.
$$
Let $K>0$ be the constant such that the Lemmas \ref{shift}, \ref{partialshift},
\ref{sizeerror}, and \ref{deformation} hold for these collections. Furthermore, choose $\delta_0>0$ such that
\begin{equation}\label{gamma}
\gamma=\left(\frac{\theta_2 \nu}{\theta_1^{1+\alpha}}\right)^4\cdot (1+K\delta_0)^{1+\alpha}<1.
\end{equation}
We may assume that for every $F\in  W^{ss}_{\text{loc}}(F_*)$ 
and  $w_1, w_2\in \{0,1\}^{4+1}$
$$
\dist(B^4_{w_1}(F), B^4_{w_2}(F))\ge \delta_0.
$$
Use the notation $d=d(F, \tilde{F})$.

\bigskip

The proof is by induction. The Lemma holds for $n=0$ because $h_0=\id$. It holds even with $C_0=0$. But we will choose  $C_0>0$  large enough such that the Lemma holds with $C=C_0$ whenever $|y-x|\ge \delta_0$. This is indeed possible because from Lemma \ref{c1conv} we get a uniform bound on the norm of $Dh_n$.

\bigskip

Assume the Lemma holds for $n$:
$$
|h_{n}(y)-h_n(x)-Dh_n(x)(y-x)|\le C_n  \cdot d \cdot |y-x|^{1+\alpha}.
$$
Choose $x,y\in \OO_{n+1}$ with $|y-x|\le \delta_0$. Then there is $w\in \{0,1\}^{4}$, and $v\in \{0,1\}$, with $x,y\in B^4_{wv}$. Let 
$$
\psi=\Psi^4_w(F)
$$
and
$$
\tilde{\psi}=\Psi^4_w(\tilde{F}).
$$
Observe,  $\psi^{-1}(x),\psi^{-1}(y)\in B_0(R^4F)\cup B_1(R^4F)$. This will allow us to apply Proposition \ref{esti}.
The maps $\Psi^4_w(F)$ depend analytically on $F$, see Theorem \ref{hyp}(8). Hence,
\begin{equation}\label{dpsipsi}
||\psi-\tilde{\psi}||_{C^2}\le C\cdot d.
\end{equation}
Let
$$
\hat{h}_n=h_n(R^4F,R^4\tilde{F})
$$
then
$$
h_{n+1}(x)=\tilde{\psi}\circ \hat{h}_n\circ \psi^{-1}(x)
$$
and
$$
h_{n+1}(y)=\tilde{\psi}\circ \hat{h}_n\circ \psi^{-1}(y).
$$
Our goal is to estimate $|h_{n+1}(y)-h_{n+1}(x)-Dh_{n+1}(x)(y-x)|$, the distortion (The appendix gives some estimates on distortions in general).
To avoid very cumbersome expressions we will work backward.

\bigskip

The maps $\psi$ and $\tilde{\psi}$ will contribute to the distortion. The first part of the proof will lead to an estimate on the influence of these distortions, see (\ref{AB}). From Lemma \ref{c0conv} we get
\begin{equation}\label{hpsi}
\hat{h}_n(\psi^{-1}(x))=\psi^{-1}(x)+\Delta
\end{equation}
with
\begin{equation}\label{del}
|\Delta|\le C\cdot d
\end{equation}
where we used $d(R^4F, R^4\tilde{F})\le \nu^4 d$, see Theorem \ref{hyp}(6). From Lemma \ref{c1conv} we get
\begin{equation}\label{Dhpsi}
D\hat{h}_n(\psi^{-1}(x))=\id +O(d).
\end{equation}
Use (\ref{hpsi}) and (\ref{Dhpsi}) and Lemma \ref{sizeerror} in the following estimates. Let
$$
\begin{aligned}
z&=\hat{h}_n(\psi^{-1}(x))+D\hat{h}_n(\psi^{-1}(x))(D\psi^{-1}(x)(y-x)+E_{\psi^{-1}}(x,y))\\
&=\psi^{-1}(x)+D\psi^{-1}(x)(y-x)+E_{\psi^{-1}}(x,y)+\Delta+O(d|y-x|).\\
\end{aligned}
$$
Hence,
\begin{equation}\label{z}
z=\psi^{-1}(y)+\Delta+O(d|y-x|).
\end{equation}
The next term we treat, using (\ref{hpsi}), (\ref{z}), Lemma \ref{deformation}, (\ref{dpsipsi}), Lemma \ref{partialshift},  is
$$
\begin{aligned}
E_{\tilde{\psi}}(\hat{h}_n(\psi^{-1}(x)),z)&=
E_{\tilde{\psi}}(\psi^{-1}(x)+\Delta, \psi^{-1}(y)+\Delta+O(d|y-x|))\\
&=E_\psi(\psi^{-1}(x)+\Delta, \psi^{-1}(y)+\Delta)+O(d|y-x|^2)\\
\end{aligned}
$$
Hence, using Lemma \ref{shift}, and (\ref{del})  we obtain
\begin{equation}\label{B}
E_{\tilde{\psi}}(\hat{h}_n(\psi^{-1}(x)),z)=E_\psi(\psi^{-1}(x), \psi^{-1}(y))+O(d|y-x|^2).
\end{equation}
The next term we treat, using (\ref{dpsipsi}), Lemma \ref{c1conv}, (\ref{hpsi}), and (\ref{del}), is
$$
\begin{aligned}
&D\tilde{\psi}(\hat{h}_n(\psi^{-1}(x)))\cdot D\hat{h}_n(\psi^{-1}(x))\cdot E_{\psi^{-1}}(x,y)=\\
&(D\psi(\hat{h}_n(\psi^{-1}(x)))+O(d))\cdot 
(\id +O(d))\cdot 
E_{\psi^{-1}}(x,y)=\\
&(D\psi(\psi^{-1}(x))+O(d))\cdot 
(\id +O(d))\cdot 
E_{\psi^{-1}}(x,y).
\end{aligned}
$$
Hence, using Lemma \ref{sizeerror},
\begin{equation}\label{A}
\begin{aligned}
&D\tilde{\psi}(\hat{h}_n(\psi^{-1}(x)))D\hat{h}_n(\psi^{-1}(x)) E_{\psi^{-1}}(x,y)=\\
&D\psi((\psi^{-1}(x))E_{\psi^{-1}}(x,y)+O(d|y-x|^2).
\end{aligned}
\end{equation}
Now we apply  the cancellation from Lemma \ref{inv} to (\ref{B}) and (\ref{A}) to obtain an estimate on the influence of the distortions of $\psi$ and $\tilde{\psi}$ 
\begin{equation}\label{AB}
D\tilde{\psi}(\hat{h}_n(x_-))D\hat{h}_n(x_-) E_{\psi^{-1}}(x,y)+
E_{\tilde{\psi}}(\hat{h}_n(x_-),z)=
O(d|y-x|^2),
\end{equation}
where $x_-=\psi^{-1}(x)$.

\bigskip

\noindent
The second part of the proof will give an estimate on the influence of the distortion of the approximate conjugation $\hat{h}_n$ between $R^4F$ and $R^4\tilde{F}$. 
First use the induction hypothesis,  to obtain
\begin{equation}\label{hp-1y}
\begin{aligned}
\hat{h}_n(\psi^{-1}(y))&=
\hat{h}_n(\psi^{-1}(x)+D\psi^{-1}(x)(y-x)+E_{\psi^{-1}}(x,y))\\
&=z+\Delta z,
\end{aligned}
\end{equation}
where
\begin{equation}\label{Dz1}
|\Delta z|\le C_n\cdot d(R^4F, R^4\tilde{F})\cdot 
| D\psi^{-1}(x)(y-x)+ E_{\psi^{-1}}(x,y)  |^{1+\alpha}.
\end{equation}
Now use Theorem \ref{hyp}(6), Proposition \ref{esti}, 
 Lemma \ref{sizeerror}, to obtain
\begin{equation}\label{Dz2}
|\Delta z|\le 
C_n\cdot \left(\frac{\nu}{\theta_1^{1+\alpha}}\right)^4\cdot d \cdot
| y-x  |^{1+\alpha}\cdot (1+K|y-x|)^{1+\alpha}.
\end{equation}
Hence, again using Proposition \ref{esti}, (\ref{gamma}), and $|y-x|\le \delta_0$,
\begin{equation}\label{dz}
|\tilde{\psi}(z+\Delta z)-\tilde{\psi}(z)|\le C_n\cdot \gamma \cdot d \cdot | y-x  |^{1+\alpha}.
\end{equation}
The final term is
$$
\begin{aligned}
\tilde{\psi}(z)=&h_{n+1}(x)+\\
  &D\tilde{\psi}(\hat{h}_n(x_-))\cdot 
D\hat{h}_n(x_-)\cdot 
\{D\psi^{-1}(x)(y-x)+E_{\psi^{-1}}(x,y)\}+\\
 &E_{\tilde{\psi}}(\hat{h}_n(\psi^{-1}(x)),z)\\
 =&h_{n+1}(x)+Dh_{n+1}(x)(y-x)+\\
 &D\tilde{\psi}(\hat{h}_n(x_-))\cdot D\hat{h}_n(x_-)\cdot E_{\psi^{-1}}(x,y)+E_{\tilde{\psi}}(\hat{h}_n(x_-),z).
 \end{aligned}
$$
From  (\ref{AB}) we obtain
\begin{equation}\label{fiz}
\tilde{\psi}(z)=h_{n+1}(x)+Dh_{n+1}(x)(y-x)+O(d|y-x|^2).
\end{equation}

\bigskip

\noindent
Finally, we will use both parts of the proof to estimate the distortion of $h_{n+1}$.
Using (\ref{hp-1y}) we can estimate 
$$
h_{n+1}(y)=\tilde{\psi}(z+\Delta z)=\tilde{\psi}(z)+(\tilde{\psi}(z+\Delta z)-\tilde{\psi}(z))
$$
by  (\ref{dz}) and (\ref{fiz}), we obtain
$$
\begin{aligned}
&|h_{n+1}(y)-h_{n+1}(x)-Dh_{n+1}(x)(y-x)|\le  \\
&C_n \cdot \gamma \cdot  d \cdot  | y-x  |^{1+\alpha}+C\cdot d\cdot |y-x|^2=\\
&(C_n \gamma+C)\cdot d \cdot | y-x  |^{1+\alpha}.
 \end{aligned}
$$
Hence,
\begin{equation}\label{contr}
C_{n+1}\le \max\{C_0, \gamma \cdot C_n +C\}.
\end{equation}
This finishes the proof of the Lemma.
\end{proof}

\begin{rem}\label{cancel dist} The application of Lemma \ref{c2conv} does not need the factor $d(F, \tilde{F})$. However, this factor plays a crucial role. The contraction in (\ref{contr}) by the factor $\gamma<1$ arises in the equation
 (\ref{Dz1}), (\ref{Dz2}), and (\ref{dz}). The distance factor reveals the contraction.
 
The proof became lengthy. The cause is the error term in (\ref{fiz}). 
This is the sum of the distortions of $\tilde{\psi}$ and $\psi^{-1}$.
The estimates show that these distortions essentially cancel each other.
Without evoking the cancelation in (\ref{AB}) the error term would be of the order 
$|y-x|^2$ without the necessary factor $d(F, \tilde{F})$. Similar analysis occurs in the proof of Lemma \ref{caDh}.
\end{rem}

\begin{lem}\label{caDh} 	For every $0<\alpha<\alpha_0$ there exists a constant $C>0$ such that 
$$
|Dh_n(y)-Dh_n(x)|\le C  \cdot d(F, \tilde{F})\cdot |y-x|^{\alpha},
$$
where $x, y\in \OO_n$.
\end{lem}

\begin{proof}  Fix $<\alpha<\alpha_0$. 
Consider the collection of diffeomorphisms 
$$
\{\Psi^4_w(F)| F\in  W^{ss}_{\text{loc}}(F_*) \text{ and } w\in \{0,1\}^4\}
$$
and
$$
\{ (\Psi^4_w(F))^{-1}| F\in  W^{ss}_{\text{loc}}(F_*) \text{ and } w\in \{0,1\}^4\}.
$$
The Lemmas \ref{shiftED}, \ref{partialshiftED},
\ref{sizeerrorED}, and \ref{deformationED} hold uniformly for these collections.  
Choose $\delta_0>0$ such that for every $F\in  W^{ss}_{\text{loc}}(F_*)$ 
and  $w_1, w_2\in \{0,1\}^{4+1}$
$$
\dist(B^4_{w_1}(F), B^4_{w_2}(F))\ge \delta_0.
$$
Use the notation $d=d(F, \tilde{F})$.

\bigskip

The proof is by induction. The Lemma holds for $n=0$ because $h_0=\id$. It holds even with $C_0=0$. But we will choose  $C_0>0$  large enough such that the Lemma holds with $C=C_0$ whenever $|y-x|\ge \delta_0$. This is indeed possible because from Lemma \ref{c1conv} we get a uniform bound on the norm of $Dh_n$.

\bigskip

Assume the Lemma holds for $n$:
$$
|Dh_{n}(y)-Dh_n(x)|\le C_n  \cdot d \cdot |y-x|^{\alpha}.
$$
Choose $x,y\in \OO_{n+1}$ with $|y-x|\le \delta_0$. Then there is $w\in \{0,1\}^{4}$, $v\in  \{0,1\}$,  with $x,y\in B^4_{wv}$. Let 
$$
\psi=\Psi^4_w(F)
$$
and
$$
\tilde{\psi}=\Psi^4_w(\tilde{F}).
$$
Observe,  $\psi^{-1}(x),\psi^{-1}(y)\in B_0(R^4F)\cup B_1(R^4F)$. This will allow us to apply Proposition \ref{esti}. Let
$$
\hat{h}_n=h_n(R^4F,R^4\tilde{F})
$$
then
$$
Dh_{n+1}(x)=D\tilde{\psi}(\hat{h}_n(\psi^{-1}(x)))
D\hat{h}_n(\psi^{-1}(x)) D\psi^{-1}(x)
$$
and
$$
Dh_{n+1}(y)=D\tilde{\psi}(\hat{h}_n(\psi^{-1}(y)))
D\hat{h}_n(\psi^{-1}(y)) D\psi^{-1}(y).
$$
Our goal is to estimate $|Dh_{n+1}(y)-Dh_{n+1}(x)|$.
We will work backward again.
The maps $\psi$ and $\tilde{\psi}$ will contribute. The first part of the proof will lead to an estimate on this influence, see (\ref{edtotal}).

\bigskip

Use (\ref{hpsi}), (\ref{Dhpsi}), Lemma \ref{c2conv}, and Proposition \ref{esti} in the following estimates. Let $x_-=\psi^{-1}(x)$ then
$$
\begin{aligned}
\hat{h}_n(\psi^{-1}(y))&=\hat{h}_n(x_-)+
D\hat{h}_n(x_-)(\psi^{-1}(y)-\psi^{-1}(x))+O(d\cdot |y-x|^{1+\alpha})\\
&=\psi^{-1}(x)+\Delta+(\psi^{-1}(y)-\psi^{-1}(x))+O(d|y-x|).\\
\end{aligned}
$$
Hence,
\begin{equation}\label{zED}
\hat{h}_n(\psi^{-1}(y))=\psi^{-1}(y)+\Delta+O(d|y-x|).
\end{equation}
Let $y_-=\psi^{-1}(y)$ then, by using (\ref{zED}), (\ref{hpsi}), Lemma \ref{deformationED}, Lemma \ref{shiftED}, and (\ref{del})
\begin{equation}\label{ED1}
ED_{\tilde{\psi}}(\hat{h}_n(x_-), \hat{h}_n(y_-))
=ED_{\psi}(x_-, y_-)+O(d|y-x|).
\end{equation}
Use (\ref{hpsi}), (\ref{del}),  Lemma \ref{sizeerrorED}
\begin{equation}\label{DED}
D\tilde{\psi}(\hat{h}_n(x_-))ED_{\psi^{-1}}(x,y)=
D\psi(x_-)ED_{\psi^{-1}}(x,y)+O(d|y-x|).
\end{equation}
The influence of $\psi$ and $\tilde{\psi}$ is estimated by using 
(\ref{ED1}) and (\ref{DED})
$$
\begin{aligned}
&D\tilde{\psi}(\hat{h}_n(x_-))ED_{\psi^{-1}}(x,y)+
ED_{\tilde{\psi}}(\hat{h}_n(x_-), \hat{h}_n(y_-))D\psi^{-1}(y)\\
&=D\psi(x_-)ED_{\psi^{-1}}(x,y)+
ED_{\psi}(x_-, y_-)D\psi^{-1}(y)+O(d|y-x|).
\end{aligned}
$$
Hence, by using the cancellation Lemma \ref{invED}, 
\begin{equation}\label{edtotal}
D\tilde{\psi}(\tilde{x}_-)ED_{\psi^{-1}}(x,y)+
ED_{\tilde{\psi}}(\tilde{x}_-, \tilde{y}_-))D\psi^{-1}(y)=
O(d|y-x|),
\end{equation}
where $\tilde{x}_-=\hat{h}_n(x_-)$
and $\tilde{y}_=\hat{h}_n(y_-)$.

\bigskip

\noindent
Finally, using Lemma \ref{c1conv}, lemma \ref{sizeerrorED}, and  (\ref{edtotal}), we are prepared to estimate
$$
\begin{aligned}
Dh_{n+1}&(y)=\\ 
=&D\tilde{\psi}(\tilde{y}_-)
D\hat{h}_n(\psi^{-1}(y)) D\psi^{-1}(y)\\
=&(D\tilde{\psi}(\tilde{x}_-)+ED_{\tilde{\psi}}(\tilde{x}_-, \tilde{y}_-))
D\hat{h}_n(y_-)(D\psi^{-1}(x)+ED_{\psi^{-1}}(x,y))\\
=&D\tilde{\psi}(\tilde{x}_-)D\hat{h}_n(y_-)D\psi^{-1}(x)+
D\tilde{\psi}(\tilde{x}_-)(\id +O(d))ED_{\psi^{-1}}(x,y)+\\
& \text{   } ED_{\tilde{\psi}}(\tilde{x}_-, \tilde{y}_-)(\id +O(d))D\psi^{-1}(y)\\
=&D\tilde{\psi}(\tilde{x}_-)D\hat{h}_n(y_-)D\psi^{-1}(x)+O(d|y-x|).
\end{aligned}
$$
Now, use the induction hypothesis, Proposition \ref{esti}, and Theorem \ref{hyp}(6) in the following estimates
$$
\begin{aligned}
|Dh_{n+1}&(y)-Dh_{n+1}(x)|\\
\le& \text{   } |D\tilde{\psi}(\tilde{x}_-)
(D\hat{h}_n(y_-)-D\hat{h}_n(x_-) )
D\psi^{-1}(x)|+
C\cdot d\cdot |y-x|\\
\le & \text{   } \theta_2^4\cdot C_n\cdot d(R^4F,R^4\tilde{F})\cdot 
\left( \frac{|y-x|}{\theta_1^4}  \right)^\alpha\cdot \frac{1}{\theta_1^4}+
C\cdot d\cdot |y-x|\\
\le & \text{    } C_n \cdot \left(\frac{\theta_2 \nu}{\theta_1^{1+\alpha}}\right)^4\cdot d\cdot |y-x|^\alpha+ C\cdot d\cdot |y-x|.
\end{aligned}
$$
Recall the definition of $\gamma<1$, see (\ref{gamma}),  and we get
$$
C_{n+1}\le \max\{C_0, \gamma \cdot C_n+C\}.
$$
The Lemma follows.
\end{proof}

\begin{rem}\label{nutheta2}
The proof of Lemma \ref{c2conv} and Lemma \ref{caDh} rely on the strong convergence in 
$W^{ss}_{\text{loc}}(F_*)$. Namely,
$$
\frac{\theta_2\nu}{\theta_1}<1.
$$
Also observe that Lemma \ref{c1conv} is only useful under the same condition. The convergence in $W^{s}_{\text{loc}}(F_*)$ is too slow 
to allow a similar treatment. 

The Cantor sets $\OO_F$ are constructed similar to the limit set of an iterated function system. See Figure \ref{micro}. The pieces $B^n_{wv}$ are images of the branches
$\Psi^n_w$. These branches are generated by the rescaling maps $\psi^n_0$ and $\psi^n_1$ which have the property $\psi^n_0\to \psi_0$ and $\psi^n_1\to \psi_1$. 
 One could generalize this construction by arbitrarily chosen converging sequences of contracting diffeomorphisms.
In this generalized context one would also obtain a rigidity theorem under the condition that the convergence is fast enough compared to the minimal and maximal contraction rates of the limit diffeomorphisms. However,
one can construct examples of this generalized type with slow convergence such that the conjugations are not Lipschitz. 

The Rigidity Theorem \ref{rigid} in the context of area-preserving maps relies on the above mentioned general principle together with the fact that the weak stable direction is a family containing the fixed point of 
analytically conjugated maps, see Theorem \ref{hyp}(7). 
\end{rem}

From Lemma \ref{c1conv} we get convergence 
$
\lim Dh_n(x)= D(x).
$
The following Lemmas states that the matrices $D(x)$ are indeed the derivatives of the limit map $h=\lim h_n:\OO_F\to \OO_{\tilde{F}}$ and that these derivatives satisfy a Holder condition. They follow from Lemma \ref{c2conv} and Lemma \ref{caDh} by taking the limit.

\begin{lem}\label{c2h} 	For every $0<\alpha<\alpha_0$ there exists a constant $C>0$ such that 
$$
|h(y)-h(x)-Dh(x)(y-x)|\le C \cdot d(F, \tilde{F})\cdot  
|y-x|^{1+\alpha},
$$
where $x, y\in \OO_F$ and $Dh(x)=D(x)$.
\end{lem}

\begin{lem}\label{HolderDh} 	For every $0<\alpha<\alpha_0$ there exists a constant $C>0$ such that 
$$
|Dh(y)-Dh(x)|\le C  \cdot d(F, \tilde{F})\cdot |y-x|^{\alpha},
$$
where $x, y\in \OO_F$.
\end{lem}

\comm{

\begin{proof} For $x,y\in \OO_F$ let
$$
E(x,y)=h(y)-h(x)-Dh(x)(y-x).
$$
From Lemma \ref{c2h} we get
$$
|E(x,y)|\le C \cdot d(F, \tilde{F})\cdot  |y-x|^{1+\alpha}.
$$
Observe,
$$
E(x,y)+E(y,x)=(Dh(y)-Dh(x))(y-x).
$$
Hence, in matrix norm, 
$$
|Dh(y)-Dh(x)|\le 2C \cdot d(F, \tilde{F})\cdot |y-x|^{\alpha} .
$$
The Lemma follows.
\end{proof}

}

\bigskip

\noindent
{\it Proof of Proposituion \ref{core}}: Lemma \ref{c2h} and Lemma \ref{HolderDh} say that the conjugation $h:\OO_F\to \OO_{\tilde{F}}$ is a $C^{1+\alpha}$-map. 
By switching the role of $F$ and $\tilde{F}$ we see that $h$ is in fact a $C^{1+\alpha}$-diffeomorphism. The Proposition follows. 
\qed

\bigskip

\noindent
{\it Proof of the Rigidity Theorem \ref{rigthm}}: Given $F\in W^{s}(F_*)$, let $\tilde{F}\in W^{ss}_{\text{loc}}(F_*)$  and $n\ge 1$ be the map from Lemma \ref{friend}. In particular, the conjugation $h:\OO_{R^nF}\to \OO_{\tilde{F}}$ is a $C^{1+\text{Lip}}$-diffeomorphism. Proposition \ref{core} says,
$$
\OO_{\tilde{F}} =\OO_{F_*} \mod(C^{1+\alpha}).
$$ Hence,
$$
\OO_{R^nF} =\OO_{F_*} \mod(C^{1+\alpha}).
$$
Apply Lemma \ref{RnFconj} and the Rigidity Theorem \ref{rigthm} follows.
\qed

%% file: appendix.tex
\section{Appendix}\label{app}

Let $D\subset \mathbb{R}^2$ and consider a bounded collection of $C^3$-maps $\psi:D\to  \mathbb{R}^2$. For such a collection there exists a $K>0$ such that the following Lemmas hold for the error terms
$$
E_\psi(x,y)=\psi(y)-\psi(x)-D\psi(x)(y-x).
$$
corresponding to affine approximations. We associate these error terms with distortion.

\begin{lem}\label{shift}  For every $x,y,x+\Delta, y+\Delta\in D$
$$
|E_{\psi}(x+\Delta,y+\Delta)-E_{\psi}(x,y)|\le K\cdot |\Delta|\cdot |y-x|^2. 
$$
\end{lem}

\begin{lem}\label{partialshift}  For every $x,y,y+\Delta\in D$
$$
|E_{\psi}(x,y+\Delta)-E_{\psi}(x,y)|\le K\cdot \{ |\Delta|\cdot |y-x|+
|\Delta|^2\}. 
$$
\end{lem}

\begin{lem}\label{sizeerror}  For every $x,y\in D$
$$
|E_{\psi}(x,y)|\le K\cdot |y-x|^2.
$$
\end{lem}

\begin{lem}\label{deformation}  For every $x,y\in D$
$$
|E_{\tilde{\psi}}(x,y)-E_{\psi}(x,y)|\le K\cdot ||\tilde{\psi}-\psi||_{C^2}\cdot |y-x|^2.
$$
\end{lem}

For completeness let us indicate how one can show these estimates. First observe that it suffices to obtain these estimates for the coordinate functions. Secondly, one can reduce the problem to a one-variable problem by restricting these coordinate function to the lines connecting two points. 
For a one-variable function $\psi:(d_1,d_2)\to \mathbb{R}$ we have 
$$
E_\psi(x,y)=\int_x^y \psi^{(2)}(t)\cdot (y-t)\cdot dt.
$$
One can use this to prove for example Lemma \ref{partialshift}. Namely,
$$
|E_{\psi}(x,y+\Delta)-E_{\psi}(x,y)|\le 
$$
$$
|  \int_y^{y+\Delta} \psi^{(2)}(t)\cdot (y-t)\cdot dt  + 
\Delta\cdot \int_x^{y+\Delta} \psi^{(2)}(t) dt |.
$$
The Lemma follows.

\bigskip

The next Lemma is responsible for a crucial cancellation in the proof of Lemma \ref{c2conv}.

\begin{lem}\label{inv} If $\psi:D\to \psi(D)$ is a diffeomorphism then 
$$
D\psi(\psi^{-1}(x))E_{\psi^{-1}}(x,y)+E_{\psi}(\psi^{-1}(x),\psi^{-1}(y))=0,
$$
for every $x,y\in \psi(D)$
\end{lem}

The proof is a calculation. Let $x_-=\psi^{-1}(x)$ and $y_-=\psi^{-1}(y)$. 
$$
\begin{aligned}
0&=\psi\circ \psi^{-1}(y)-y\\
&= \psi(\psi^{-1}(x)+D\psi^{-1}(x)(y-x)+E_{\psi^{-1}}(x,y))-y\\
&=x+D\psi(x_-) \{D\psi^{-1}(x)(y-x)+E_{\psi^{-1}}(x,y)\}+
E_{\psi}(x_-,y_-)-y\\
&=D\psi(\psi^{-1}(x))E_{\psi^{-1}}(x,y)+E_{\psi}(\psi^{-1}(x),\psi^{-1}(y)).
\end{aligned}
$$

\bigskip

Consider the error terms
$$
ED_\psi(x,y)=D\psi(y)-D\psi(x).
$$

\begin{lem}\label{shiftED}  For every $x,y,x+\Delta, y+\Delta\in D$
$$
|ED_{\psi}(x+\Delta,y+\Delta)-ED_{\psi}(x,y)|\le K\cdot |\Delta|\cdot |y-x|. 
$$
\end{lem}

\begin{lem}\label{partialshiftED}  For every $x,y,y+\Delta\in D$
$$
|ED_{\psi}(x,y+\Delta)-ED_{\psi}(x,y)|\le K\cdot  |\Delta|. 
$$
\end{lem}

\begin{lem}\label{sizeerrorED}  For every $x,y\in D$
$$
|ED_{\psi}(x,y)|\le K\cdot |y-x|.
$$
\end{lem}

\begin{lem}\label{deformationED}  For every $x,y\in D$
$$
|ED_{\tilde{\psi}}(x,y)-ED_{\psi}(x,y)|\le K\cdot ||\tilde{\psi}-\psi||_{C^2}\cdot |y-x|.
$$
\end{lem}

The next Lemma is responsible for a crucial cancellation in the proof of Lemma \ref{caDh}.

\begin{lem}\label{invED} If $\psi:D\to \psi(D)$ is a diffeomorphism then 
$$
D\psi(\psi^{-1}(x))ED_{\psi^{-1}}(x,y)+ED_{\psi}(\psi^{-1}(x),\psi^{-1}(y))D\psi^{-1}(y)=0,
$$
for every $x,y\in \psi(D)$
\end{lem}

%% file: references.tex
\comm{
\harvarditem{Abad {\it et al}}{2000}{AK} J. J. Abad, H. Koch, 
Renormalization and periodic orbits for Hamiltonian flows, {\it Comm.  
Math.  Phys. } { \bf 212} (2000) \# 2 371--394.

\harvarditem{Abad  {\it et al}}{1998}{AKW} J. J. Abad, H. Koch and P. 
Wittwer,  A renormalization group for Hamiltonians: numerical results, 
{\it Nonlinearity} { \bf 11} (1998) 1185--1194.

\harvarditem{Benettin  {\it et al}} {1980} {BCGG} G. Benettin et al, 
Universal properties in conservative dynamical systems, {\it Lettere al 
Nuovo Cimento} {\bf 28} (1980) 1--4.

\harvarditem{Bountis}{1981}{Bou} T. Bountis,  Period doubling 
bifurcations and universality in conservative Systems, {\it Physica} { 
\bf 3D} (1981) 577--589.

\harvarditem {CAPD} {2009} {CAPD} CAPD--Computer Assisted Proofs in Dynamics group, a C++ package for rigorous numerics,
{\tt http://capd.wsb-nlu.edu.pl}

\harvarditem{de Carvalho {\it et al}}{2005}{dCLM} A. de Carvalho, M. 
Lyubich, M. Martens,  Renormalization in the H\'enon family, I: 
Universality but non-rigidity, {\it J. Stat. Phys} { \bf 121} (2005) 
611--669.

\harvarditem{Collet {\it et al}}{1980}{CEK1} P. Collet, J.-P. Eckmann 
and H. Koch, Period doubling bifurcations for families of maps on 
${\fR}^n$, {\it J. Stat. Phys.} { \bf 3D} (1980).

\harvarditem{Collet  {\it et al}}{1981}{CEK2} P. Collet, J.-P. Eckmann 
and H. Koch, On universality for area-preserving maps of the plane , 
{\it Physica} { \bf 3D} (1981) 457--467.

\harvarditem{Derrida and Pomeau}{1980}{DP} B. Derrida, Y. Pomeau, 
Feigenbaum's ratios of two dimensional  area preserving maps, {\it Phys. 
Lett.} {\bf A80} (1980) 217--219.

\harvarditem{Duarte}{2000}{Duarte1} P. Duarte, Persistent homoclinic 
tangencies for conservative maps near identity, {\it Ergod. Th. \& 
Dynam. Sys.} {\bf 20} (2000) 393--438.

\harvarditem{Eckmann  {\it et al}}{1982}{EKW1} J.-P. Eckmann, H. Koch 
and P. Wittwer, Existence of a fixed point of the doubling 
transformation for area-preserving maps of the plane, {\it Phys. Rev. A} 
{ \bf 26} (1982) \# 1  720--722.

\harvarditem{Eckmann  {\it et al}}{1984}{EKW2} J.-P. Eckmann, H. Koch 
and P. Wittwer, A Computer-Assisted Proof of Universality for 
Area-Preserving Maps, {\it Memoirs of the American Mathematical Society} 
{\bf  47} (1984), 1--121.

\harvarditem{Epstein}{1986}{Eps1} H. Epstein, New proofs of the 
existence of the Feigenbaum  functions, {\it  Commun. Math. Phys.} {\bf 
106} (1986) 395--426.

\harvarditem{Epstein}{1989}{Eps2} H. Epstein, Fixed points of 
composition operators II, {\it  Nonlinearity } {\bf 2} (1989) 305--310.

\harvarditem{Escande and Doveil}{1981}{ED} D. F. Escande, F. Doveil, 
Renormalization method for computing the threshold of the large scale 
stochastic instability in two degree of freedom Hamiltonian systems, 
{\it  J .Stat. Phys. } {\bf 26} (1981) 257--284. 

\harvarditem{de Faria}{1992}{dF1} E. de Faria, {\it Proof of 
universality for critical circle mappings}, Thesis, CUNY, 1992.

\harvarditem{de Faria}{1999}{dF2} E. de Faria, Asymptotic rigidity of 
scaling ratios for critical circle mappings, {\it Ergodic Theory Dynam. 
Systems} {\bf 19} (1999), no. 4, 995--1035.

\harvarditem{Feigenbaum}{1978}{Fei1} M. J. Feigenbaum,  Quantitative 
universality for a class of nonlinear transformations, {\it J. Stat. 
Phys.} {\bf 19} (1978) 25--52.

\harvarditem{Feigenbaum}{1979}{Fei2} M. J. Feigenbaum, Universal metric 
properties of non-linear transformations, {\it J. Stat. Phys.} {\bf 21} 
 (1979) 669--706.

\harvarditem{Gaidashev}{2005}{Gai1} D. Gaidashev, Renormalization of 
isoenergetically degenerate Hamiltonian flows and associated 
bifurcations of invariant tori, {\it Discrete Contin. Dyn. Syst.} {\bf 
13} (2005), no. 1, 63--102.

\harvarditem{Gaidashev}{2007}{Gai3}  D. Gaidashev,  Cylinder 
renormalization for Siegel disks and a constructive Measurable Riemann 
Mapping Theorem, {\it  Nonlinearity} {\bf 20} (2007), no 3, 713--742.

\harvarditem{Gaidashev}{2010}{Gai4}  D. Gaidashev,  Period Doubling Renormalization for Area-Preserving Maps and Mild Computer Assistance in Contraction Mapping Principle,
e-print {\tt math.DS/1009.0625} at {\tt Arxiv.org}.

\harvarditem{Gaidashev and Johnson}{2009a}{GJ1}D. Gaidashev, T. 
Johnson, Dynamics of the Universal Area-Preserving Map Associated with Period Doubling: Hyperbolic Sets, \textit{Nonlinearity} {\bf 22} 2487-2520.

\harvarditem{Gaidashev and Johnson}{2009b}{GJ2}D. Gaidashev, T. Johnson, Dynamics of the Universal Area-Preserving Map Associated with Period Doubling: Stable Sets, \textit{J. Mod. Dyn.} {\bf 3} (2009), no 4, 555--587.

\harvarditem{Gaidashev and Koch}{2004}{GK} D. Gaidashev, H. Koch, 
Renormalization and shearless invariant tori: numerical results, {\it 
Nonlinearity} {\bf 17} (2004), no. 5, 1713--1722.

\harvarditem{Gaidashev and Koch}{2008}{GK1} D. Gaidashev, H. Koch, 
Period doubling in area-preserving maps: an associated one-dimenisonal 
problem, e-print {\tt math.DS/0811.2588} at {\tt Arxiv.org}.

\harvarditem{Gaidashev and Yampolsky}{2007}{GaiYa}D. Gaidashev, M. 
Yampolsky, Cylinder renormalization of Siegel disks, {\it Exp. Math.} 
{\bf 16:2} (2007).

%\harvarditem {Gorodetsky and Kaloshin} {2008} {GoKa} A. Gorodetsky and 
%V. Kaloshin, Conservative homoclinic bifurcations and some applications, 
%preprint (2008).

\harvarditem{Helleman}{1980}{Hel} R. H. G. Helleman, Self-generated 
chaotic behavior in nonlinear mechanics, in "Fundamental problems in 
statistical mechanics", Ed. by E. G. D. Cohen, North-Holland, Amsterdam, 
p.165, (1980).

\harvarditem{Johnson}{2010}{J1} T. Johnson, No elliptic islands for the universal area-preserving map, 
e-print {\tt math.DS/1012.2409} at {\tt Arxiv.org}.

\harvarditem{Katok and Hasselblatt}{1995}{KH} A. Katok, B. Hasselblat, 
Introduction to the Modern Theory of Dynamical Systems, Cambridge University 
Press, Cambridge (1995).

\harvarditem{Khanin {\it et al}}{2007}{KLDM}  K. Khanin, J. Lopes Dias, 
J. Marklof,  Multidimensional continued fractions, dynamic 
renormalization and KAM theory,  {\it Comm. Math. Phys.},  {\bf 270}  
 (2007),  no. 1, 197--231.

\harvarditem{Koch}{2002}{Koch1}H. Koch, On the renormalization of 
Hamiltonian flows, and critical invariant tori, {\it Discrete Contin.  
Dyn.  Syst.}  {\bf 8} (2002),  633--646.

\harvarditem{Koch}{2004}{Koch2}H. Koch, A renormalization group fixed 
point associated with the breakup of golden invariant tori, {\it  
Discrete Contin. Dyn. Syst.}  {\bf 11}  (2004),  no. 4, 881--909.

\harvarditem{Koch}{2008}{Koch3} H. Koch, Existence of critical 
invariant tori, \textit{Ergodic theory and dynamical systems} {\bf 28} 1879-94

\harvarditem{Koch {\it et al}}{1996}{KWS} H. Koch, A. Schenkel, P. Wittwer,  Computer-Assisted Proofs in Analysis and Programming in Logic: A Case Study, \textit{SIAM Review} {\bf 38}(4) (1996) 

\harvarditem{Koci\'c}{2005}{Kocic} S. Koci\'c, Renormalization of 
Hamiltonians for Diophantine frequency vectors and KAM tori,{\it 
Nonlinearity} {\bf 18} (2005) 2513--2544.

\harvarditem{Kokubu {\it et al}} {2007} {KWZ07} H. Kokubu, D. Wilczak and P. Zgliczy\'nski, Rigorous verification of cocoon bifurcations in the Michelson system,  {\it Nonlinearity}  {\bf 20}  no. 9 (2007), 2147--2174.

\harvarditem{Lyubich}{1999}{Lyu} M. Lyubich, Feigenbaum-Coullet-Tresser 
universality and Milnor's  hairness conjecture, {\it Annals of 
Mathematics} {\bf 149} (1999) 319--420.

\harvarditem{Lyubich and Martens}{2008}{LM} M. Lyubich, M. Martens,  
Renormalization in the H\'enon family, II: Homoclinic tangle, to appear in {\it Inventiones Mathematicae}. preprint, arXiv:mathDS/0804.0780
(2008).

\harvarditem{MacKay}{1982}{McK1} R. S. MacKay, Renormalisation in area 
preserving maps, Thesis, Prin\-ce\-ton (1982). World Scientific, London 
 (1993).

\harvarditem{MacKay}{1983}{McK2} R. S. MacKay, Renormalization approach 
to invariant circles in area-preserving maps, {\it Physica} {\bf D7} 
 (1983) 283--300.

\harvarditem{Mehr and Escande}{1984}{ME} A. Mehr and D.F. Escande, 
Destruction of KAM tori in Hamiltonian systems: link with the 
distabilization of nearby cycles and calculation of residues, {\it  
Physica }{\bf D13}  (1984) 302--338.

\harvarditem{de Melo and Pinto}{1999}{dMP} W. de Melo, A. A. Pinto, Rigidity of $C^2$ infinitely renormalizable unimodal maps, {\it Comm. Math. Phys},  {\bf 208} (1999), 91--105. 

\harvarditem{McMullen}{1996}{McM1} C. McMullen, {\it Renormalization and 3-manifolds which fiber over the circle. (Annals of Mathematics Studies, 142)}, Princeton University Press, Princeton, NJ, 1996.

\harvarditem{McMullen}{1998}{McM} C. McMullen, Self-similarity of Siegel 
disks and Hausdorff dimension of  Julia sets, {\it Acta Math.} {\bf 
180} (1998), 247--292.

\harvarditem{Neumaier}{1990}{Ne90} A.\,Neumaier,
Interval Methods for Systems of Equations.
Encyclopedia of Mathematics and its Applications 37,
Cambridge Univ. Press, Cambridge, 1990

\harvarditem{Pa\l uba}{1989}{Pal} W. Pa\l uba, The Lipschitz condition for the conjugacies of Feigenbaum-like mappings, {\it  Fund. Math.}, {\bf  132} (1989), 227--258. 

%\harvarditem{Palis and Takens}{1993}{PT} J. Palis, F. Takens, {\it  
%Hyperbolicity and sensitive chaotic dynamics at homoclinic bifurcations. 
%Fractal dimensions and infinitely many attractors. (Cambridge Studies in 
%Advanced Mathematics, 35.)} Cambridge University Press, Cambridge, 1993.

\harvarditem{Rand}{1988}{Rand}  D. A. Rand, Global phase space universality, smooth conjugacies and renormalisation. I. The $C^{1+\alpha}$ case, {\it Nonlinearity}, {\bf   1} (1988), 181--202. 

\harvarditem{Shenker and Kadanoff}{1982}{Shen} S. J. Shenker, L. P. 
Kadanoff, Critical behaviour of KAM surfaces. I Empirical results, {\it  
J. Stat.  Phys. } {\bf 27} (1982) 631--656.

\harvarditem{Sullivan}{1992}{Sul} D. Sullivan, Bounds, quadratic 
differentials and renormalization conjectures, in: Mathematics into the 
Twenty-first Century, AMS Centennial Publications, Vol. II, Amer. Math. 
Soc., Providence, R.I. (1992) 417-466.

\harvarditem{Tresser and Coullet}{1978}{TC} C. Tresser and P. Coullet, 
It\'erations d'endomorphismes et groupe de renormalisation, {\it  C. R. 
Acad.  Sci. Paris} {\bf 287A} (1978), 577--580.

\harvarditem{Yampolsky}{2002}{Ya1} M. Yampolsky, Hyperbolicity of 
renormalization of critical circle maps,
{\it  Publ. Math. Inst. Hautes Etudes Sci.} {\bf 96} (2002), 1--41.

\harvarditem{Yampolsky}{2003}{Ya2} M. Yampolsky, Renormalization 
horseshoe for critical circle maps,
{\it  Commun. Math. Physics} {\bf 240} (2003), 75--96.

\harvarditem{Yampolsky}{2007}{Ya3} M. Yampolsky,  Siegel disks and 
renormalization fixed points, \textit{Holomorphic Dynamics and Renormalization (Fields Inst. Commun.} vol 53)
(Providence, RI:AMS) pp. 377-93.

\harvarditem {Zgliczy\'nski} {1997} {Z97} P. Zgliczy\'nski, Computer 
assisted proof of chaos in the R\"ossler equations and in the H\'enon 
map,  {\it Nonlinearity}  {\bf 10},  no. 1 (1997), 243--252.

\harvarditem {Zgliczy\'nski} {2009} {Z09} P. Zgliczy\'nski, Covering 
relations, cone conditions and the stable manifold theorem, {\it Journal 
of Differential Equations} {\bf 246} issue 5 (2009), 1774--1819.

\harvarditem {Zgliczy\'nski and Gidea} {2004} {ZG04} P. Zgliczy\'nski 
and M. Gidea 2004 Covering relations for multidimensional dynamical 
systems.  J. Differential Equations  {\bf 202}  (2004),  no. 1, 32--58.

\harvarditem {Progs} {2009} {GP} Programs available at {\tt http://math.uu.se/$\thicksim$gaidash}
\harvarditem {Progs} {2011} {JP} Programs available at {\tt http://www.math.cornell.edu/$\thicksim$tjohnson}

\end{thebibliography}